\documentclass{amsart}
\usepackage{amsfonts}
\usepackage{amsmath,amssymb}
\usepackage{amsthm}
\usepackage{amscd}
\usepackage{graphics}
\usepackage{graphicx}

\theoremstyle{remark}{
\newtheorem{Def}{{\rm Definition}}
\newtheorem{Ex}{{\rm Example}}
\newtheorem{Rem}{{\rm Remark}}

}
\theoremstyle{plain}
{

\newtheorem{Prop}{Proposition}

\newtheorem{MainThm}{Main Theorem}

}

\begin{document}
\title[On real algebraic maps topologically special generic maps]{Some remark on real algebraic maps which are topologically special generic maps and generalize the canonical projections of the unit spheres}
\author{Naoki kitazawa}

\address{Institute of Mathematics for Industry, Kyushu University, 744 Motooka, Nishi-ku Fukuoka 819-0395, Japan\\
 TEL (Office): +81-92-802-4402 \\
 FAX (Office): +81-92-802-4405 \\
}
\email{n-kitazawa@imi.kyushu-u.ac.jp, naokikitazawa.formath@gmail.com}
\urladdr{https://naokikitazawa.github.io/NaokiKitazawa.html}
\maketitle
\begin{abstract}

Morse functions with exactly two singular points on spheres and projections of naturally embedded spheres in the Euclidean spaces (the {\it canonical projections of the unit spheres}) are generalized to {\it special generic} maps. This is a smooth map which is the composition of a surjection onto a manifold with a smooth immersion of codimension $0$. The preimage of the surjection is diffeomorphic to the unit sphere of a suitable dimension at each point of the interior of the manifold and a single point set at each point on the boundary.

We construct real algebraic maps topologically seen as special generic maps whose images are smoothly embedded manifolds. The canonical projections of the unit spheres are generalized here. We also constructed some real algebraic maps topologically special generic maps previously. We construct such maps in more general situations and give additional remarks.  We are also interested in explicit and meaningful smooth maps in differential topology and real algebraic geometry. The construction is important and difficult.

\end{abstract}
\section{Introduction.}
\label{sec:1}

In theory of Morse functions, related singularity theory and applications to geometry of manifolds, Morse functions with exactly two singular points on spheres (in Reeb's theorem) are of most fundamental functions and maps. Such maps are generalized to {\it special generic} maps: a special generic map $c$ is a smooth map locally represented as a projection or by the form $c(x_1,\cdots x_m)=(x_1,\cdots x_{n-1},{\Sigma}_{j=1}^{m-n+1} {x_j}^2)$.
 \cite{burletderham, furuyaporto, saeki2} are related pioneering studies. For Morse functions, see \cite{reeb} for a pioneering study and \cite{milnor2} for systematic exposition on fundamental notions and arguments.

This paper also studies real algebraic geometry. Complex algebraic geometry is presenting one of central topics in algebraic geometry.
Real algebraic geometry has another nice history.
\cite{nash, tognoli} show classical and important studies and \cite{bochnakcosteroy, bochnakkucharz, kollar, kucharz} present fundamental notions and arguments and advanced studies, for example. Our paper studies explicit construction in real algebraic geometry, respecting and generalizing (the canonical projections of the unit) spheres.

Celebrated theory by Nash, followed by Tognoli, says that smooth closed manifolds are non-singular real algebraic manifolds. They are also shown to be the zero sets of some real polynomial maps. 
We can also approximate smooth maps between {\it non-singular} real algebraic manifolds by real algebraic ones. 
It is also important and difficult to have real algebraic functions and maps explicitly. 

For example, \cite{maciasvirgospereirasaez, ramanujam, takeuchi} show explicit cases and some fundamental properties of real algebraic functions in the theory of Lie groups and symmetric spaces. The unit spheres are simplest symmetric spaces. We consider another generalization for the projections of spheres: constructing real algebraic maps topologically special generic maps. The author has succeeded in such construction in \cite{kitazawa2}, followed by \cite{kitazawa5} by the author. We present new related result. More precisely, we construct new explicit maps on some real algebraic hypersurfaces. As one of related future work, \cite{gabard, viro} present closely related long-standing studies on construction of real algebraic curves in the plane and general real algebraic hypersurfaces (Remark \ref{rem:6}).
\subsection{Smooth manifolds and maps.}
Let $X$ be a topological space with the structure of a cell complex whose dimension is finite. We can define the dimension $\dim X$ uniquely, which is a non-negative integer. 
A topological manifold is well-known to be homeomorphic to some CW complex. A smooth manifold is well-known to be a polyhedron canonically, and a so-called PL manifold. It is also well-known that a topological space homeomorphic to a polyhedron whose dimension is at most $2$ has the structure of a polyhedron uniquely. If we restrict this to the class of topological manifolds, then the condition is weaken to be that the dimension is at most $3$. Consult \cite{moise} for example. Hereafter, for a manifold $X$, let $\partial X$ denote the boundary of it and ${\rm Int}\ X:=X -\partial X$ the interior of it.
  
Let ${\mathbb{R}}^k$ denote the $k$-dimensional Euclidean space, which is one of most fundamental $k$-dimensional smooth manifolds. It is also the Riemannian manifold with the standard Euclidean metric. 
Let $\mathbb{R}:={\mathbb{R}}^1$. 
For $x \in {\mathbb{R}}^k$, $||x|| \geq 0$ is defined as the canonically defined distance between $x$ and the origin $0$.
This is also regarded as a ({\it non-singular}) real algebraic manifold or the $k$-dimensional {\it real affine space}. $S^k:=\{x \in {\mathbb{R}}^{k+1} \mid ||x||=1\}$ denotes the $k$-dimensional unit sphere, which is a $k$-dimensional smooth compact submanifold of ${\mathbb{R}}^{k+1}$ and has no boundary. It is connected for any positive integer $k \geq 1$. It is a discrete set with exactly two points for $k=0$. It is the zero set of the real polynomial $||x||^2-1={\Sigma}_{j=1}^{k+1} {x_j}^2-1$ ($x:=(x_1,\cdots,x_{k+1})$) and a (non-singular) real algebraic (sub)manifold (hypersurface). 
$D^k:=\{x \in {\mathbb{R}}^{k} \mid ||x|| \leq 1\}$ is the $k$-dimensional unit disk. It is a $k$-dimensional smooth compact and connected submanifold of ${\mathbb{R}}^{k}$. Let $c:X \rightarrow Y$ be a differentiable map between differentiable manifolds $X$ and $Y$. $x \in X$ is a {\it singular} point of $c$ if the rank of the differential at $x$ is smaller than the minimum of $\{\dim X, \dim Y\}$ and $c(x)$ is called a {\it singular value} of $c$ for a singular point $x$ of $c$. Let $S(c)$ denote the set of all singular points of $c$. Here, most of differentiable maps are smooth, equivalently, of the class $C^{\infty}$.
The canonical projection of the Euclidean space ${\mathbb{R}}^k$ is denoted
by ${\pi}_{k,k_1}:{\mathbb{R}}^{k} \rightarrow {\mathbb{R}}^{k_1}$. This maps 
each point $x=(x_1,x_2) \in {\mathbb{R}}^{k_1} \times {\mathbb{R}}^{k_2}={\mathbb{R}}^k$ to $x_1 \in {\mathbb{R}}^{k_1}$ with $k_1, k_2>0$ and $k=k_1+k_2$. The canonical projection of the unit sphere $S^{k-1}$ is its restriction, rigorously.

Here, a {\it real algebraic} manifold is the union of some connected components of the zero sets of a real polynomial map. {\it Non-singular} real algebraic manifolds are defined naturally by the implicit function theorem for the polynomial maps. The real affine space and the unit sphere are of simplest examples. Our {\it real algebraic} maps are the compositions of the canonical embeddings into the real affine spaces with the canonical projections. Note that they may be so-called {\it Nash} manifolds and maps, as in \cite{nash}. However, we do not use such terminologies.
\subsection{Our main result and the structure of our paper.}


Spaces of smooth maps between two manifolds are topologized with so-called {\it $C^{\infty}$ Whitney topologies}. Roughly, two smooth maps are close if their differentials of arbitrary degrees are close. {\it $C^r$ Whitney topologies} are weaker topologies, respecting differentials of degrees at most $r$. The $C^{r_1}$ Whitney topologies are weaker than the $C^{r_2}$ Whitney topologies for $0 \leq r_1<r_2$: $r_1$ is a non-negative integer and $r_2$ is a positive integer or $\infty$. These topologies are important in singularity theory and differential topology. Rigorous exposition on singularity theory of differential maps including these topologies is left to \cite{golubitskyguillemin}. {\it Diffeomorphisms} are smooth homeomorphisms with no singular points (of them). The {\it diffeomorphism group} of a smooth manifold is the group of all diffeomorphisms from the manifold onto itself. A {\it smooth} bundle is a bundle whose fiber is a smooth manifold and whose structure group is (a subgroup of) the diffeomorphism group. A {\it linear} bundle is a smooth bundle whose fiber is the Euclidean space, the unit disk, or the unit sphere and whose structure group consists of linear transformations, defined naturally. 

Proposition \ref{prop:1} explains a fundamental theorem on structures of special generic maps. Using the notation, we present our new real algebraic map which is special generic topologically in Main Theorem \ref{mthm:1}, proven in the next section with our additional result Main Theorems \ref{mthm:2} and \ref{mthm:3}. Main Theorem \ref{mthm:2} is for additional comments on Main Theorem \ref{mthm:1}. Main Theorem \ref{mthm:3} generalizes them. The third section presents additional comments on our new result.  

\begin{Prop}[\cite{saeki2}]
	\label{prop:1}
	Let $m>n \geq 1$ be integers. For a special generic map $f:M \rightarrow N$ from an $m$-dimensional closed manifold $M$ into an $n$-dimensional smooth manifold $N$ with no boundary, we have a decomposition $f=\bar{f} \circ q_f$ as follows. Note that the case $N:={\mathbb{R}}^n$ is mainly discussed in \cite{saeki2} and our paper, and that as is also presented shortly in \cite{saeki2}, $N$ can be general.  
	\begin{enumerate}
		\item There exist an $n$-dimensional smooth compact manifold $W_f$, a smooth surjection $q_f:M \rightarrow W_f$, mapping $S(f)$ onto $\partial W_f$ by a diffeomorphism, and a smooth immersion $\bar{f}$. 
		\item There exist a small collar neighborhood $N(\partial W_f)$ of $\partial W_f$ and a linear bundle whose projection is the composition of the restriction of $q_f$ to the preimage of $N(\partial W_f)$ with the canonical projection to $\partial W_f$ and whose fiber is $D^{m-n+1}$. By {\rm (}the restriction of{\rm )} $q_f$, we also have a smooth bundle over the closure of $W_f-N(\partial W_f)$ in $W_f$ whose fiber is the unit sphere $S^{m-n}$.
	\end{enumerate}
\end{Prop}

\begin{MainThm}
	\label{mthm:1}
	Let $m \geq n \geq 1$ be integers. Suppose that a smooth embedding ${\bar{f}}_{{\mathbb{R}}^n}:\bar{N} \rightarrow {\mathbb{R}}^n$ from an $n$-dimensional compact and connected manifold $\bar{N}$ into ${\mathbb{R}}^n$ exists. Then we have an $m$-dimensional non-singular real algebraic manifold $M$ being also compact, connected, and the zero set of a real polynomial function and the following maps.
	\begin{enumerate}
		\item A special generic map $f:M \rightarrow N:={\mathbb{R}}^n$ in Proposition \ref{prop:1} with $\bar{f}={\bar{f}}_{{\mathbb{R}}^n}:W_f=\bar{N} \rightarrow {\mathbb{R}}^n$. This map $f$ is not real algebraic in general.
		\item A real algebraic map $f_0:M \rightarrow {\mathbb{R}}^n$ and a homeomorphism $\phi:{\mathbb{R}}^n \rightarrow {\mathbb{R}}^n$ satisfying the relation $f_0=\phi \circ f$ and $\phi(f(S(f)))=f_0(S(f_0))${\rm :} remember that $S(c)$ denotes the singular set of a smooth map $c$.
	\end{enumerate}
\end{MainThm}
 The simplest case $\bar{N}=D^n$ accounts for the canonical projection $f=f_0:S^m \rightarrow {\mathbb{R}}^n$ of the unit sphere $S^m$ with $\phi$ being the identity map.

\noindent {\bf Conflict of interest.} \\
This work was also supported by the Research Institute of Mathematical Sciences, an International Joint Usage/Research Center located in Kyoto University: the conference "Extension of the singularity theory" (https://www2.akita-nct.ac.jp/kasedou/workshop/rims2023/index.html)".
The author works at Institute of Mathematics for Industry (https://www.jgmi.kyushu-u.ac.jp/en/about/young-mentors/). This is closely related to our study. Our study thanks this project for supports. The author is also a researcher at Osaka Central
Advanced Mathematical Institute (OCAMI researcher), supported by MEXT Promotion of Distinctive Joint Research Center Program JPMXP0723833165. He is not employed there. However this helps our studies. Our study also thanks this. \\
\ \\
{\bf Data availability.} \\
Data essentially supporting our present study are all here. This paper is regarded as an improved version of one of a
preprints of the author \cite{kitazawa6}.
\section{On Main Theorems.}
For $a_1,a_2 \in \mathbb{R}$ such that $a_1<a_2$, let the closed interval denoted by $[a_1,a_2]:=\{x \in \mathbb{R} \mid a_1 \leq x \leq a_2\}$ as is the standard notation. 
Let the open interval denoted by $(a_1,a_2):={\rm Int}\ [a_1,a_2]$ as is the standard case. Similarly, we also define $(a_1,a_2]:=\{x \in \mathbb{R} \mid a_1<x \leq a_2\}$, for example and for $a \in \mathbb{R}$, let us define $[a,\infty):=\{x \in \mathbb{R} \mid x \geq a\}$ for example.

As an important ingredient in special generic maps, we explain the {\it Morse function for the natural height of the disk $D^k \subset {\mathbb{R}}^k$}. We map $x \in D^k$ to $(\sqrt{1-{||x||}^2},x) \in \mathbb{R} \times D^k \subset {\mathbb{R}}^{k+1}$ by a diffeomorphism, compose the projection to the second component and we have the {\it Morse function for the natural height of the disk $D^k \subset {\mathbb{R}}^k$}. 
The projection of $D^k$ embedded in ${\mathbb{R}}^{k+1}$ here is the restriction of ${\pi}_{k+1,1} {\mid}_{S^k}$ to $S^k \bigcap \{x:=(x_1,\cdots x_{k+1}) \mid x_1 \geq 0\}$.
Locally, a special generic map from an $m$-dimensional manifold with no boundary into an $n$-dimensional manifold with no boundary with $m \geq n$ is a smooth map such that around each singular point of it, the map is represented as the product map of such a function and the identity map on an ($n-1$)-dimensional manifold (where suitable local coordinates are chosen). 

Related to Proposition \ref{prop:1}, we also review a theorem reconstructing special generic maps from given smooth immersions of codimension $0$, conversely. Note again that the case $N:={\mathbb{R}}^n$ is discussed in the original study mainly, and that we can consider the general case.
In \cite{kitazawa7}, followed by \cite{kitazawa4}, the author has redefined several notions and consider specific reconstruction. We also abuse these notions for example.

\begin{Prop}[\cite{saeki2}]
\label{prop:2}
Let $m \geq n \geq 1$ be integers. Suppose that a smooth immersion ${\bar{f}}_N:\bar{N} \rightarrow N$ from an $n$-dimensional compact manifold $\bar{N}$ into an $n$-dimensional smooth manifold $N$ with no boundary and the following additional data {\rm (}objects{\rm )} are given. Hereafter, we can abuse the notation in Proposition \ref{prop:1} and do. We can also abuse this in the case $m=n$ and we do.

\begin{itemize}
\item A small collar neighborhood $N(\partial \bar{N})$ of $\partial \bar{N}$ and a smooth bundle ${\bar{M}}_{\rm Interior}={\bar{M}}_{\rm I}$ over the closure $\overline{\bar{N}- N(\partial \bar{N})}$ of $\bar{N}- N(\partial \bar{N})$ in $\bar{N}$ whose fiber is the unit sphere $S^{m-n} \subset \partial D^{m-n+1}$. We call this bundle an {\rm interior unit sphere bundle for ${\bar{f}}_N$}.
\item A linear bundle ${\bar{M}}_{\rm Boundary}={\bar{M}}_{\rm B}$ over $\partial \bar{N}$ whose fiber is the unit disk $D^{m-n+1}$. We call this bundle a {\rm boundary linear bundle for ${\bar{f}}_N$}.
\item We consider the smooth bundle $\partial {\bar{M}}_{\rm I}$ as the restriction to the boundary $\partial (\overline{\bar{N}- N(\partial \bar{N})})$ of the closure $\overline{\bar{N}- N(\partial \bar{N})}$. We also consider the subbundle $\partial {\bar{M}}_{\rm B}$ whose fiber is $S^{m-n}=\partial D^{m-n+1}$. We consider the identification between their base spaces in the natural way. More precisely, we can regard the collar neighborhood $N(\partial \bar{N})$ of $\partial \bar{N}$ as the product $(\partial \bar{N}) \times [0,1]$. Furthermore, each point $x \in \partial \bar{N}$ is naturally identified with a point $(x,0) \in (\partial \bar{N}) \times \{0\}$, of course. In addition, in considering the base space of ${\bar{M}}_{\rm Boundary}={\bar{M}}_{\rm B}$ over $\partial N(\partial \bar{N})-\partial \bar{N}=(\partial \bar{N}) \times \{1\}$ as $\partial \bar{N}$, $x \in \partial \bar{N}$ and $(x,0) \in (\partial \bar{N}) \times \{0\}$ are also seen as $(x,1) \in (\partial \bar{N}) \times \{1\}$. 

In addition, the manifold $\partial N(\partial \bar{N})-\partial \bar{N}=(\partial \bar{N}) \times \{1\}$ is also regarded as the subspace of $\overline{\bar{N}- N(\partial \bar{N})} \subset \bar{N}$ and these two spaces are naturally identified. The third object is a bundle isomorphism between the smooth bundle $\partial {\bar{M}}_{\rm I}$ over $\partial ( \overline{\bar{N}- N(\partial \bar{N})}) \subset  \overline{\bar{N}- N(\partial \bar{N})}$ and the smooth bundle $\partial {\bar{M}}_{\rm B}$ over $\partial \bar{N}=\partial N(\partial \bar{N})-\partial \bar{N}$ where the induced map between the base spaces is the natural identification for the spaces of course{\rm :} note also that we identify several spaces according to the identification before. We call this bundle isomorphism an {\rm identification bundle isomorphism for ${\bar{f}}_N$}.

\end{itemize}
Then there exist an $m$-dimensional closed manifold $M$
and a special generic map $f:M \rightarrow N$ such that $\bar{f}$ in Proposition \ref{prop:1} and ${\bar{f}}_N$ agree. We can also have the manifold and the map enjoying the following properties.
\begin{enumerate}
\item \label{prop:2.1}
The interior unit sphere bundle ${\bar{M}}_{\rm I}$ for ${\bar{f}}_N$ and the bundle over $\overline{\bar{N}-N(\partial \bar{N})}=\overline{W_f-N(\partial W_f)}$ in Proposition \ref{prop:1} are isomorphic and naturally identified via a diffeomorphism ${\Phi}_{M,{\rm I}}:{\bar{M}}_{\rm I} \rightarrow M-{q_f}^{-1}(N(\partial \bar{N}))$ and the identity map on $\overline{\bar{N}-N(\partial \bar{N})}=\overline{W_f-N(\partial W_f)}$ where the closure is taken in $\bar{N}=W_f$ as before.
\item \label{prop:2.2}
The linear bundle whose fiber is $D^{m-n+1}$ in Proposition \ref{prop:1} and the boundary linear bundle ${\bar{M}}_{\rm B}$ for ${\bar{f}}_N$ are isomorphic and naturally identified via a diffeomorphism ${\Phi}_{M,{\rm B}}:{\bar{M}}_{\rm B} \rightarrow {q_f}^{-1}(N(\partial \bar{N}))$ and the identity map on $N(\partial \bar{N})=N(\partial W_f)$.
\item \label{prop:2.3}
We consider the bundle over $\overline{\bar{N}-N(\partial \bar{N})}=\overline{W_f-N(\partial W_f)}$ and the linear bundle whose fiber is $D^{m-n+1}$ in Proposition \ref{prop:1}. Their restrictions to the boundaries are identified by the identification induced canonically from the manifold $M$ and the map $f:M \rightarrow N$.

We also consider the restriction of the smooth bundle ${\bar{M}}_{\rm I}$ in {\rm (}\ref{prop:2.1}{\rm )} to the boundary and the subbundle of the linear bundle ${\bar{M}}_{\rm B}$ whose fiber is $\partial D^{m-n+1}$ in {\rm (}\ref{prop:2.2}{\rm )}. We consider their identification or the bundle isomorphism induced from {\rm (}\ref{prop:2.1}, \ref{prop:2.2}{\rm )} and the short arguments on $M$ and $f:M \rightarrow N$ just before{\rm :} the identification coincides with the identification bundle isomorphism for ${\bar{f}}_N$.

\end{enumerate}
\end{Prop}
\begin{Rem}
	\label{rem:1}
	This remarks reconstruction of a special generic map $f:M \rightarrow N$ in Proposition \ref{prop:2}.
	\begin{enumerate}
		\item \label{rem:1.1} 
By considering the identification bundle isomorphism for ${\bar{f}}_N$ and gluing (the spaces of the domains of) ${\Phi}_{M,{\rm I}}$ and ${\Phi}_{M,{\rm B}}$ by this, we have a diffeomorphism ${\Phi}_{M,{\rm I}, {\rm B}}$ onto our new manifold $M$ such that ${\Phi}_{M,{\rm I}}$ and ${\Phi}_{M,{\rm B}}$ are regarded as its restrictions.
\item \label{rem:1.2} We can regard that the interior unit sphere bundle ${\bar{M}}_{\rm I}$ for ${\bar{f}}_N$ and the boundary linear bundle ${\bar{M}}_{\rm B}$ for ${\bar{f}}_N$ naturally reconstructs a smooth surjection $q_f:M \rightarrow W_f=\bar{N}$. After that, we compose the given smooth immersion ${\bar{f}}_N:\bar{N} \rightarrow N$ to have our desired special generic map $f:M \rightarrow N$.
\end{enumerate} 
\end{Rem}

Note that in Proposition \ref{prop:1}, we cannot consider the condition "$m \geq n$". The case with the constraint $m=n$ is closely related to Eliashberg's celebrated theory \cite{eliashberg1}. See also \cite{eliashberg2}.

\begin{Def}
	\label{def:1}
In Proposition \ref{prop:2}, if the interior unit sphere bundle  for ${\bar{f}}_N$ and the boundary linear bundle for ${\bar{f}}_N$ are trivial and the identification bundle isomorphism for ${\bar{f}}_N$ is the product map of the identity map on the base space (or the identification between the base spaces) and the identity map on the fiber $S^{m-1}=\partial D^{m-n+1}$, then it is called a {\it product-organized} special generic map.
\end{Def}

\cite{kitazawa7} has first introduced such a class of special generic maps. We first name the class in this way in \cite{kitazawa4}.
Of course in our present paper, we do not need to understand precise arguments or results from \cite{kitazawa4, kitazawa7} well.

We note the canonical projection of the unit sphere is of such a class.
We consider a manifold represented as a connected sum $\sharp (S^{k_j} \times S^{m-k_j})$ of finitely many manifolds $S^{k_j} \times S^{m-k_j}$ where $m \geq n \geq 2$ and $1 \leq k_j \leq n-1$ are integers. We take the connected sum in the smooth category. We can reconstruct a product-organized special generic map into ${\mathbb{R}}^n$ there from a natural embedding ${\bar{f}}_{{\mathbb{R}}^n}$ whose image is diffeomorphic to the boundary connected sum $\natural (D^{k_j} \times S^{n-k_j})$ of finitely many manifolds $D^{k_j} \times S^{n-k_j}$ taken in the smooth category. This gives another simplest example.

We have the following by the definition and the structures of the maps immediately. In the following, we abuse the notation of Propositions \ref{prop:1} and \ref{prop:2} or its natural variant.

\begin{Prop}
	\label{prop:3}
	Let $f_1:M \rightarrow N$ be a product-organized special generic map and for the induced smooth immersion $\bar{f}:=\bar{f_1}:W_{f_1} \rightarrow N$ and another smooth immersion ${{\bar{f}}_N}^{\prime}:W_{f_1} \rightarrow N$, let there exist a homeomorphism ${\phi}_{N_{1,2}}:N \rightarrow N$ satisfying the relation ${{\bar{f}}_N}^{\prime}={\phi}_{N_{1,2}} \circ \bar{f_1}$. In this situation, there exists a product-organized special generic map $f_2:M \rightarrow N$ satisfying the relation $f_2= {\phi}_{N_{1,2}} \circ f_1$ and $\bar{f_2}={{\bar{f}}_N}^{\prime}$. 
	\end{Prop}
In short, Remark \ref{rem:1} (\ref{rem:1.2}) gives Proposition \ref{prop:3}. We explain some, more precisely.
First, from $f_1$, we have two kinds of bundles as in Proposition \ref{prop:1}. From these two bundles and the identification between the boundaries, we obtain the three objects for reconstruction of a special generic map in Proposition \ref{prop:2} canonically. We can consider a smooth surjection as in Remark \ref{rem:1} (\ref{rem:1.2}), and after that, we compose the immersion ${{\bar{f}}_N}^{\prime}$ to have our desired special generic map $f_2$.

The following new result is a kind of our remark on Main Theorem \ref{mthm:1}.

\begin{MainThm}
\label{mthm:2}
In Main Theorem \ref{mthm:1}, $f$ and $f_0$ can be obtained as
product-organized special generic maps.
\end{MainThm}

\begin{proof}[A proof of Main Theorems \ref{mthm:1} and \ref{mthm:2}]

STEPs 1--3 are main ingredients. \\
\ \\
STEP 1 Obtaining another nice compact manifold which is isotopic to the original embedding and whose boundary consists of connected components each of which is the zero set of some real polynomial function and non-singular. \\

This is a new ingredient in our related studies. Our arguments here thanks to \cite{kollar}. Especially, "Special Case 5" and "Discussion 14". 

By considering a suitable small isotopy of the class $C^1$ starting from the original embedding ${\bar{f}}_{{\mathbb{R}}^n}$, we have another nice smooth embedding ${\bar{f}}_{{\mathbb{R}}^n,0}:\bar{N} \rightarrow {\mathbb{R}}^n$ in the following. Before explaining the new embedding, we give some small remark. We can also consider the isotopy in the class $C^a$ where $a$ is an arbitrary integer greater than $1$ according to \cite{lellis}. For this approximation, see also \cite{bodinpopescupampusorea, elredge}. This can be also chosen as an isotopy fixing the manifold except a small collar neighborhood $N(\partial \bar{N})$ of the boundary $\partial \bar{N}$.

Let ${\bar{N}}_{0} \subset {\mathbb{R}}^n$ denote the resulting image ${\bar{f}}_{{\mathbb{R}}^n,0}(\bar{N})$. Let ${\partial}_j {\bar{N}}_{0}$ denote all connected components of the boundary $\partial {\bar{N}}_{0}$, indexed by all positive integers $j=1,\cdots,l \geq 1$ smaller than or equal to $l$. The manifold ${\partial}_j {\bar{N}}_{0}$ is a non-singular real algebraic manifold of dimension $n-1$ and the zero set of a real polynomial $f_{{\partial}_j}$.
In addition, we can also get the embedding and the subset with the following properties.
\begin{itemize}
\item ${\bar{N}}_{0} = {\bigcap}_{j=1}^l \{x \in {\mathbb{R}}^n \mid f_{{\partial}_j}(x) \geq 0\}$.
\item ${\rm Int}\ {\bar{N}}_{0} = {\bigcap}_{j=1}^l \{x \in {\mathbb{R}}^n \mid f_{{\partial}_j}(x) > 0\}$.
\item At each point $x$ of ${\mathbb{R}}^n-{\bar{N}}_{0}$, the relation $f_{{\partial}_{j_x}}(x)< 0$ holds for exactly one integer $1 \leq j_x \leq l$ and the relation $f_{{\partial}_j}(x)> 0$ holds for the remaining polynomials.
\end{itemize}

\ \\
STEP 2 Arguments from \cite{kitazawa2, kitazawa6} to have our smooth real algebraic map $f_0$. \\

We apply arguments from main ingredients of \cite{kitazawa2, kitazawa6}. \\
We consider the zero set $S_0:=\{(x,y)=(x,y_1,\cdots y_{m-n+1}) \in {\mathbb{R}}^n \times {\mathbb{R}}^{m-n+1}={\mathbb{R}}^{m+1} \mid {\prod}_{j=1}^l (f_{{\partial}_j}(x))-{\Sigma}_{j=1}^{m-n+1} {y_j}^2=0\}$.

The last part of STEP 1 implies the fact that $S_0=\{(x,y) \in {\bar{N}}_{0} \times {\mathbb{R}}^{m-n+1} \subset {\mathbb{R}}^{m+1} \mid {\prod}_{j=1}^l (f_{{\partial}_j}(x))-{\Sigma}_{j=1}^{m-n+1} {y_j}^2=0\}$.
We apply the implicit function theorem to see that this set is a smooth, closed and connected manifold and also a non-singular real algebraic manifold of dimension $m$.

Around a point of the form $(x,y) \in ({\rm Int}\ {\bar{N}}_{0}) \times {\mathbb{R}}^{m-n+1}$ in this subset, it is represented as the graph of a canonically obtained smooth function where some $y_j$ is regarded as the variable of the target of the function.

Around a point of the form $(x,y) \in (\partial {\bar{N}}_{0}) \times {\mathbb{R}}^{m-n+1}$ in this subset, we consider the implicit function theorem for the polynomial ${\prod}_{j=1}^l (f_{{\partial}_j}(x))-{\Sigma}_{j=1}^{m-n+1} {y_j}^2$ in a different way.

By considering the partial derivative for some variable $x_{j_0}$ in $x=(x_1,\cdots x_n)$, we have a value which is not $0$. We discuss this more precisely.
If $x \in {\partial}_{i_0} {\bar{N}}_{0}$ holds, then the partial derivative of $f_{{\partial}_{i_0}}(x)$ is not $0$ for some variable $x_{j_0}$, since each ${\partial}_i {\bar{N}}_{0}$ is obtained by some suitable approximation of an original connected component being also a smooth compact manifold with no boundary. Remember that this is also obtained via an isotopy in the $C^1$ class. For the variable $x_{j_0}$, the partial derivative of the polynomial ${\prod}_{j=1}^l (f_{{\partial}_j}(x))-{\Sigma}_{j=1}^{m-n+1} {y_j}^2$ is the value represented as the product of the partial derivative of $f_{{\partial}_{i_0}}(x)$ (for $x_{j_0}$) and the product of the values of the remaining polynomials $f_{{\partial}_j}(x)$. This is thanks to elementary formula in calculus and the property of the zero sets of our polynomial functions. 
Remember that $f_{{\partial}_j}(x)$ is not $0$ for any $j$ except $j=i_0$.
The value of the partial derivative of the polynomial ${\prod}_{j=1}^l (f_{{\partial}_j}(x))-{\Sigma}_{j=1}^{m-n+1} {y_j}^2$ (for $x_{j_0}$) is not $0$. 

We have checked that the set $S_0$ is a smooth, closed and connected manifold and also a non-singular real algebraic manifold of dimension $m$. Put $M_0:=S_0$.
Our desired map $f_0$ is defined as the composition of the canonical embedding of $M_0 \subset {\mathbb{R}}^{m+1}$ with the canonical projection ${\pi}_{m+1,n}:{\mathbb{R}}^{m+1} \rightarrow {\mathbb{R}}^n$. \\
\ \\
STEP 3 Is $f_0$ special generic? \\
We can have a small collar neighborhood $N(\partial {\bar{N}}_0)$ of the boundary $\partial {\bar{N}}_0$, regarded as $(\partial {\bar{N}}_0) \times [0,1] \subset N(\partial {\bar{N}}_0)$. More precisely, we can have one and consider our identification in the following way thanks to arguments related to \cite[Discussion 14]{kollar}.
\begin{itemize}
\item The manifold $(\partial {\bar{N}}_0) \times \{0\}$ is identified with $(\partial N(\partial {\bar{N}}_0)) \bigcap \partial {\bar{N}}_0 = \partial {\bar{N}}_0$.
\item The manifold $(\partial {\bar{N}}_0) \times \{t_0\} \subset (\partial {\bar{N}}_0) \times (0,1]$ is identified with the zero set of the polynomial ${\prod}_{j=1}^l (f_{{\partial}_j})-a_0 t_0$ for a sufficiently small number $a_0>0$. We have ${\Sigma}_{j=1}^{m-n+1} {y_j}^2=at_0={\prod}_{j=1}^l (f_{{\partial}_j}(x))$ at each point $(x,y) \in (((\partial {\bar{N}}_0) \times \{t_0\}) \times {\mathbb{R}}^{m-n+1}) \bigcap M_0 \subset ({\bar{N}}_0 \times {\mathbb{R}}^{m-n+1}) \bigcap M_0$ ($y=(y_1,\cdots,y_{m-n+1})$).
\end{itemize}

This means that around $N(\partial {\bar{N}}_0)$, $f_0$ is a special generic map and represented as the product map of the Morse function for the natural height of the ($m-n+1$)-dimensional unit disk $D^{m-n+1}$ and the identity map on $\partial {\bar{N}}_0$ (for suitable local coordinates).

Around a point of the form $(x,y) \in (({\rm Int}\ {\bar{N}}_{0}) \times {\mathbb{R}}^{m-n+1}) \bigcap M_0$, remember that locally $M_0$ is represented as the graph of a canonically obtained smooth function. Some $y_j$ is regarded as the variable of the target of the function here.


We can map each point of the form $(x,y) \in ({\rm Int}\ {\bar{N}}_{0}) \times ({\mathbb{R}}^{m-n+1}-\{0\})$ to $(x,\frac{1}{\sqrt{{\prod}_{j=1}^l (f_{{\partial}_j}(x))}} y)$. 
This gives a diffeomorphism on $({\rm Int}\ {\bar{N}}_{0}) \times ({\mathbb{R}}^{m-n+1}-\{0\})$. We map each point of the form $(x,0) \in ({\rm Int}\ {\bar{N}}_{0}) \times {\mathbb{R}}^{m-n+1}$ to $(x,0)$. Thus we have a diffeomorphism from $({\rm Int}\ {\bar{N}}_{0}) \times {\mathbb{R}}^{m-n+1}$ onto $({\rm Int}\ {\bar{N}}_{0}) \times {\mathbb{R}}^{m-n+1}$. By the restriction, we also have another diffeomorphism $(({\rm Int}\ {\bar{N}}_{0}) \times {\mathbb{R}}^{m-n+1}) \bigcap M_0$ onto $({\rm Int}\ {\bar{N}}_{0}) \times S^{m-n} \subset ({\rm Int}\ {\bar{N}}_{0}) \times {\mathbb{R}}^{m-n+1}$.

The structures of the manifolds and maps and these arguments imply that $f_0$ is a product-organized special generic map.  \\
\ \\
We argue STEPs 1--3 again in a generalized way in the proof of Main Theorem \ref{mthm:3}.  \\
\ \\
We check the existence of desired maps for Main Theorems \ref{mthm:1} and \ref{mthm:2}.
If we do not need to deform the embedding ${\bar{f}}_{{\mathbb{R}}^n}$ in STEP 1 and we can regard ${\bar{f}}_{{\mathbb{R}}^n}={\bar{f}}_{{\mathbb{R}}^n,0}$, then we can have desired maps by putting the special generic map $f:=f_0$ and the homeomorphism $\phi$ as the identity map. 

In general, we need to deform the embedding in STEP 1 and in the general case, our desired map $f$ is reconstructed by Propositions \ref{prop:2} and \ref{prop:3} and Definition \ref{def:1} similarly from the original embedding ${\bar{f}}_{{\mathbb{R}}^n}$. For the original embedding, each connected component of the boundary of the image ${\bar{f}}_{{\mathbb{R}}^n}(\bar{N})$ may not be real algebraic (as a result $f$ may not be real algebraic). We can also know the existence of the homeomorphism $\phi$ for Main Theorems \ref{mthm:1} and \ref{mthm:2} by the structures of the manifolds and maps and Proposition \ref{prop:3}, easily.

This completes the proof.

\end{proof}

We show a kind of generalizations for Main Theorems \ref{mthm:1} and \ref{mthm:2}.
Before this, we review several known facts on topologies and differentiable structures of spheres and disks shortly. See \cite{milnor1} and see also \cite{kervairemilnor, smale1, smale2}, for example.
A smooth manifold which is homeomorphic to $S^k$ ($k \geq 1$) is diffeomorphic to $S^k$ for $k=1,2,3,5,6$ for example. It is well-known that smooth manifolds homeomorphic to $S^k$ are classified via abstract and sophisticated theory for $k \geq 7$.
If a smooth manifold which is homeomorphic to $S^k$ is smoothly embedded into ${\mathbb{R}}^{k+1}$ for $k \geq 1$ with $k \neq 4$, then this is diffeomorphic to $S^k$. 
A smooth manifold "which is homeomorphic to $S^4$ and not diffeomorphic to $S^4$" is still undiscovered. 
A smooth manifold which is homeomorphic to $D^k$ and whose boundary is diffeomorphic to $\partial D^k=S^{k-1}$ is known to be diffeomorphic to $D^k$ unless $k=4$.
A smooth manifold which is homeomorphic to $D^k$ is known to be diffeomorphic to $D^k$ unless $k=4,5$.
Remember that "${\Delta}^k$ which is homeomorphic to $D^k$ and is not diffeomorphic to $D^k$" is still undiscovered. If this holds for some ${\Delta}^k$, then of course $k$ must be $4$ or $5$.

\begin{MainThm}
\label{mthm:3}
Suppose that a situation same as one of Main Theorem \ref{mthm:1} is given.

Let $F_{S^{m-n}}:{\mathbb{R}}^{m-n+2}={\mathbb{R}}^{m-n+1} \times \mathbb{R} \rightarrow \mathbb{R}$ be a function defined canonically by a real polynomial enjoying the following.
We naturally define $F_{S^{m-n},t}:{\mathbb{R}}^{m-n+1} \rightarrow \mathbb{R}$ as the restriction to ${\mathbb{R}}^{m-n+1} \times \{t\}$ for $t \in \mathbb{R}$ where ${\mathbb{R}}^{m-n+1} \times \{t\}$ and ${\mathbb{R}}^{m-n+1}$ are identified by the map mapping $(y,t)$ to $y$.
\begin{enumerate}
\item
\label{mthm:3.1}
 $F_{S^{m-n}}(y,t)=F_{0,S^{m-n}}(t)-F_{S^{m-n},{\mathbb{R}}^{m-n+1}}(y)$ for some real polynomials $F_{0,S^{m-n}}$ and $F_{S^{m-n},{\mathbb{R}}^{m-n+1}}$ and $(y,t) \in {\mathbb{R}}^{m-n+1} \times \mathbb{R}$.
\item
\label{mthm:3.2}
 $F_{0,S^{m-n}}(0)=0$. $F_{0,S^{m-n}}(t)>0$ for $t>0$ and $F_{0,S^{m-n}}(t)<0$ for $t<0$. The value of the derivative of $F_{0,S^{m-n}}$ at $0$ is not zero.
\item
\label{mthm:3.3}
 ${F_{S^{m-n},{\mathbb{R}}^{m-n+1}}}^{-1}(0) =\{0\} \subset {\mathbb{R}}^{m-n+1}$.
$F_{S^{m-n},{\mathbb{R}}^{m-n+1}}(y) \geq 0$ for  $y \in {\mathbb{R}}^{m-n+1}$.

\item
\label{mthm:3.4}
A positive number $a>0$ and a smooth manifold ${\Sigma}^{m-n}$ with a homeomorphism ${\phi}_{\Sigma}:{\Sigma}^{m-n} \rightarrow S^{m-n}$ are given.
 There exists a homeomorphism ${\Phi}_{F_{{S}^{m-n}},a}:{\Sigma}^{m-n} \times (0,a] \rightarrow {\bigcup}_{t \in  (0,a]} {F_{S^{m-n},t}}^{-1}(0) \subset {\mathbb{R}}^{m-n+1}$ mapping each sphere ${\Sigma}^{m-n} \times \{t\}$ onto ${F_{{S}^{m-n},t}}^{-1}(0)$ by a homeomorphism. Furthermore, the zero set ${F_{S^{m-n},t}}^{-1}(0) \subset {\mathbb{R}}^{m-n+1}$ of the real polynomial function is non-singular for $t \in (0,a]$.


\end{enumerate}

In this situation, by defining $$M_0:=\{(x,y) \in {\mathbb{R}}^n \times {\mathbb{R}}^{m-n+1}={\mathbb{R}}^{m+1} \mid F_{S^{m-n}}(y,\frac{{\prod}_{j=1}^l (f_{{\partial}_j}(x))}{T})=0\}$$
for a sufficiently large $T>0$ and real polynomials $f_{{\partial}_j}$ in our proof of Main Theorems \ref{mthm:1} and \ref{mthm:2}, we have an $m$-dimensional smooth manifold homeomorphic to $M_0$ and the following maps, similar to those in Main Theorem \ref{mthm:1}, as follows.

\begin{enumerate}
	\setcounter{enumi}{4}
	\item \label{mthm:3.5}
	 A special generic map $f:M \rightarrow N:={\mathbb{R}}^n$ in Proposition \ref{prop:1} with $\bar{f}={\bar{f}}_{{\mathbb{R}}^n}:W_f=\bar{N} \rightarrow {\mathbb{R}}^n$.
	\item \label{mthm:3.6}
	A real algebraic map $f_0:M_0 \rightarrow {\mathbb{R}}^n$ satisfying the relation $f_0 \circ \Phi=\phi \circ f $ for some pair $(\Phi:M \rightarrow M_0,\phi:{\mathbb{R}}^n \rightarrow {\mathbb{R}}^n)$ of homeomorphisms with the relations $\Phi(S(f))=S(f_0)$ and $\phi(f(S(f)))=f_0(S(f_0))$.
\end{enumerate}

Explicitly, we can also have the map $f_0:M_0 \rightarrow {\mathbb{R}}^n$ with the following properties.
\begin{enumerate}
\setcounter{enumi}{6}
\item \label{mthm:3.7}

We can abuse the notation in Proposition \ref{prop:1} replacing "$f$" by "$f_0$" for example. We have a smooth bundle whose projection is the composition of the restriction of $q_{f_0}$ to the preimage of a suitable small collar neighborhood $N(\partial W_{f_0})$ with the canonical projection to $\partial W_{f_0}$, whose fiber is homeomorphic to $D^{m-n+1}$, and the boundary of whose fiber is homeomorphic to $S^{m-n}$. Furthermore, the smooth bundle is a trivial smooth bundle. 

\item \label{mthm:3.8} The non-singular real algebraic manifold $M_0$ bounds a compact and bounded submanifold $W_0$ of dimension $m+1$ in ${\mathbb{R}}^{m+1}$ homeomorphic to a manifold obtained in the following steps.
Here ${\Delta}^k$ is for a $k$-dimensional smooth manifold homeomorphic to $D^k$ {\rm (}$k \geq 1${\rm )}. This is diffeomorphic to $D^k$ for $k \neq 4,5$. This is also diffeomorphic to $D^k$ for $k \neq 4$ if its boundary is diffeomorphic to $S^{k-1}$.
\begin{enumerate}
\item Prepare the product $\bar{N} \times {\Delta}^{m-n+1}$ and the product $(\partial \bar{N}) \times {\Delta}^{m-n+2}$. 
\item Attach the previous two manifold by a bundle isomorphism between the trivial smooth bundle $(\partial \bar{N}) \times {\Delta}^{m-n+1}$ over $\partial \bar{N}$ and the trivial smooth bundle $(\partial \bar{N}) \times {\Delta}^{m-n+1} \subset (\partial \bar{N}) \times \partial {\Delta}^{m-n+2}$ over $\partial \bar{N}$ represented as the product of the identity maps on $\partial \bar{N}$ and ${\Delta}^{m-n+1}${\rm :} here ${\Delta}^{m-n+1}$ is smoothly embedded in $\partial {\Delta}^{m-n+2}$.
\item Eliminate the corner.
\end{enumerate}
\end{enumerate}
\end{MainThm}

Hereafter, for a real number $t$ and $x \in {\mathbb{R}}^k$, we use the notation respecting the structure of the canonical vector space for ${\mathbb{R}}^k$ such as $tx \in {\mathbb{R}^k}$. The vector $tx$ is also a point of ${\mathbb{R}}^k$ of course. 

\begin{proof}[A proof of Main Theorem \ref{mthm:3}]

We can argue as "STEP 1 in the proof of Main Theorems \ref{mthm:1} and \ref{mthm:2}". We can abuse the notation there and we do.

We choose a sufficiently large number $T>0$ such that the maximal value of $\frac{{\prod}_{j=1}^l (f_{{\partial}_j}(x))}{T}$ is smaller than the given positive number $a>0$.
We define $M_0$ as presented in the statement of Main Theorem \ref{mthm:3}.

We also define the map $f_0:M_0 \rightarrow {\mathbb{R}}^n$ as the composition of the canonical embedding to ${\mathbb{R}}^{m+1}$ with the canonical projection ${\pi}_{m+1,n}$ as in Main Theorems \ref{mthm:1} and \ref{mthm:2} and their proof.

We present arguments similar to ones of "STEP 2 in the proof of Main Theorems \ref{mthm:1} and \ref{mthm:2}".
Here we also respect the assumptions (\ref{mthm:3.1})--(\ref{mthm:3.4}) here. 

First we also have the relation $$M_0:=\{(x,y) \in \bar{N_0} \times {\mathbb{R}}^{m-n+1} \mid F_{S^{m-n}}(y,\frac{{\prod}_{j=1}^l (f_{{\partial}_j}(x))}{T})=0\}$$ by our definitions and the assumptions (\ref{mthm:3.1})--(\ref{mthm:3.3}) similarly. 
The corresponding argument "in the proof of Main Theorems \ref{mthm:1} and \ref{mthm:2}" is naturally generalized.

We check that $M_0$ is also a smooth closed and connected manifold and a non-singular real algebraic manifold of dimension $m$.

Around a point of the form $(x,y) \in (({\rm Int}\ {\bar{N}}_{0}) \times {\mathbb{R}}^{m-n+1}) \bigcap M_0$, it is represented as the graph of a canonically obtained smooth function where some $y_j$ is regarded as the variable of the target of the function. We discuss this more precisely. We consider the partial derivative of the polynomial $F_{S^{m-n}}(y,\frac{{\prod}_{j=1}^l (f_{{\partial}_j}(x))}{T})=F_{0,S^{m-n}}(\frac{{\prod}_{j=1}^l (f_{{\partial}_j}(x))}{T})-F_{S^{m-n},{\mathbb{R}}^{m-n+1}}(y)$ for some variable $y_j$. According to the property (\ref{mthm:3.4}), ${F_{S^{m-n},t}}^{-1}(0) \subset {\mathbb{R}}^{m-n+1}$ is the zero set of the real polynomial which is non-singular for $t \in (0,a]$. 
We have checked that the partial derivative of the polynomial $F_{S^{m-n}}(y,\frac{{\prod}_{j=1}^l (f_{{\partial}_j}(x))}{T})=F_{0,S^{m-n}}(\frac{{\prod}_{j=1}^l (f_{{\partial}_j}(x))}{T})-F_{S^{m-n},{\mathbb{R}}^{m-n+1}}(y)$ is not $0$ for some variable $y_j$.

Around a point of the form $(x,y) \in ((\partial {\bar{N}}_{0}) \times {\mathbb{R}}^{m-n+1}) \bigcap M_0$, we have a partial derivative as in "STEP 2 in the proof of Main Theorems \ref{mthm:1} and \ref{mthm:2}". The value is the product of the derivative of $F_{0,S^{m-n}}$ at $0$, which is not $0$ from the assumption (\ref{mthm:3.2}), $\frac{1}{T}$, and the "value in the case of Main Theorems \ref{mthm:1} and \ref{mthm:2}". 
The "value in the case of Main Theorems \ref{mthm:1} and \ref{mthm:2}" means the value of the partial derivative of ${\prod}_{j=1}^l (f_{{\partial}_j}(x))$ for the "variable $x_{j_0}$ there”.
The resulting value is not $0$.

We have checked similarly that $M_0$ is a smooth closed and connected manifold and a non-singular real algebraic manifold of dimension $m$.

We present an argument similar to one in "STEP 3 in the proof of Main Theorems \ref{mthm:1} and \ref{mthm:2}".
For the boundary of the image ${\bar{N}}_{0}$, we consider a small collar neighborhood $N(\partial {\bar{N}}_{0})$. For example, we can consider this as the collar neighborhood of "STEP 3 in the proof of Main Theorems \ref{mthm:1} and \ref{mthm:2}". We regard this as $(\partial {\bar{N}}_{0}) \times [0,1]$ where $\partial {\bar{N}}_{0}$ and $(\partial {\bar{N}}_{0}) \times \{0\}$ are identified by the correspondence identifying $x$ with $(x,0)$.

We construct a homeomorphism ${\Phi}^{\prime}$ from the manifold $M_0$ onto another new smooth manifold $M$. 

First we map $(\partial {\bar{N}}_{0}) \times \{0\}$ to $(\partial {\bar{N}}_{0}) \times \{0\} \subset {\mathbb{R}}^n \times {\mathbb{R}}^{m-n+1}$. We use the natural identification mapping $(x,0)$ to $(x,0)$ here. 
According to the assumption (\ref{mthm:3.3}), we have $((\partial {\bar{N}}_{0}) \times {\mathbb{R}}^{m-n+1}) \bigcap M_0=(\partial {\bar{N}}_{0}) \times \{0\}$.


Hereafter, let ${\rm Pr}_{{\Sigma}^{m-n},a}:{\Sigma}^{m-n} \times (0,a] \rightarrow {\Sigma}^{m-n}$ denote the projection to the first component.

We can consider a smooth function $g:(0,\infty) \rightarrow (0,1]$ enjoying the following. 

\begin{itemize}
\item $g(t_1) \leq g(t_2)$ for $t_1<t_2$.
\item $g(t)=\sqrt{1-{(1-t)}^2}$ for a sufficiently small $t>0$.
\item $g(t)=1$ for $t \geq 1-\epsilon$ with a sufficiently small $\epsilon>0$.
\end{itemize}

Remember that our original homeomorphism ${\Phi}_{F_{S^{m-n}},a}$ is presented in the assumption (\ref{mthm:3.4}): this maps ${\Sigma}^{m-n} \times (0,a]$ onto ${\bigcup}_{t \in  (0,a]} {F_{S^{m-n},t}}^{-1}(0) \subset {\mathbb{R}}^{m-n+1}$ mapping each homotopy sphere ${\Sigma}^{m-n} \times \{t\}$ onto ${F_{{S}^{m-n},t}}^{-1}(0)$ by a homeomorphism.
We can map each point of the form $(x,y) \in (({\rm Int}\ {\bar{N}}_{0}) \times {\mathbb{R}}^{m-n+1}) \bigcap M_0$ to $$(x,g(t) {\phi}_{\Sigma} \circ {\rm Pr}_{{\Sigma}^{m-n},a} \circ {{\Phi}_{F_{S^{m-n}},a}}^{-1}(y))$$ according to the following rule.


\begin{itemize}
\item If $x \in N(\partial {\bar{N}}_{0})-\partial N(\partial {\bar{N}}_{0})=(\partial {\bar{N}}_{0}) \times (0,1)$, then $t$ is the value of the projection to the 2nd component of $x$ here.
\item If $x \in {\bar{N}}_{0}-{\rm Int}\ N(\partial {\bar{N}}_{0})$, then $g(t)=1$.
\end{itemize}

Thus we have a continuous bijection ${\Phi}^{\prime}$ between compact and connected manifolds, which is also a homeomorphism ${\Phi}^{\prime}$ from the smooth manifold $M_0$ onto another smooth manifold $M$. 

We can also have the composition of the canonical embedding $M \subset {\mathbb{R}}^{m+1}$ with the canonical projection ${\pi}_{m+1,n}:{\mathbb{R}}^{m+1} \rightarrow {\mathbb{R}}^n$.
The resulting map is also special generic by the construction. It is also product-organized. We define a special generic map $f^{\prime}:M \rightarrow {\mathbb{R}}^n$ as this map. Note that this is not real algebraic (in general).

An essentially same argument also appears in Main Theorems \ref{mthm:1} and \ref{mthm:2} and their proof. We further present another argument, an argument similar to which also appears there.

We define a new set $W_0 \subset {\mathbb{R}}^{m+1}$ as $$\{(x,y)=(x,y_1,\cdots y_{m-n+1}) \in {\mathbb{R}}^n \times {\mathbb{R}}^{m-n+1}={\mathbb{R}}^{m+1} \mid t \leq \frac{{\prod}_{j=1}^l (f_{{\partial}_j}(x))}{T}, F_{S^{m-n}}(y,t)=0\}$$ and by our assumption, this can be also defined as \\ $\{(x,y) \in {\bar{N}}_{0} \times {\mathbb{R}}^{m-n+1} \subset {\mathbb{R}}^{m+1} \mid 0 \leq t \leq \frac{{\prod}_{j=1}^l (f_{{\partial}_j}(x))}{T}, F_{S^{m-n}}(y,t)=0\}$.

We use the function $g$ before. We can map each point of the form $(x,y) \in ({\rm Int}\ ({\bar{N}}_{0}) \times ({\mathbb{R}}^{m-n+1}-\{0\})) \bigcap W_0$ to $(x,g(t){\phi}_{\Sigma} \circ {\rm Pr}_{{\Sigma}^{m-n},a} \circ {{\Phi}_{F_{S^{m-n}},a}}^{-1}(y))$ according to the following rule. 

\begin{itemize}
\item If $x \in N(\partial {\bar{N}}_{0})-\partial N(\partial {\bar{N}}_{0})=(\partial {\bar{N}}_{0}) \times (0,1)$, then $t$ is the value of the projection to the 2nd component of $x$ here.
\item If $x \in {\bar{N}}_{0}-{\rm Int}\ N(\partial {\bar{N}}_{0})$, then $g(t)=1$.
\end{itemize}
This rule is in fact same as the previous one.
This gives a continuous bijection. We map each point of the form $(x,0) \in {\bar{N}}_{0} \times {\mathbb{R}}^{m-n+1}$ to $(x,0)$. Thus we have a continuous bijection from $W_0$ to another compact and connected manifold, extending the homeomorphism ${\Phi}^{\prime}$ to a new homeomorphism. The image of the homeomorphism is the disjoint union of the image of the homeomorphism just before and ${\bar{N}}_{0} \times \{0\} \subset {\bar{N}}_{0} \times {\mathbb{R}}^{m-n+1}$. For here, remember again the condition (\ref{mthm:3.3}), assuming ${F_{S^{m-n},{\mathbb{R}}^{m-n+1}}}^{-1}(0) =\{0\} \subset {\mathbb{R}}^{m-n+1}$.

From this, we can see that the set $W_0-((\partial {\bar{N}}_{0}) \times \{0\})$ is regarded as the total space of a trivial bundle over ${\rm Int}\ {\bar{N}}_{0}$ which may not be a smooth bundle and whose fiber is homeomorphic to $D^{m-n+1}$. The projection is given by the restriction of the canonical projection ${\pi}_{m+1,n}$. In other words, it is also regarded as a "{\it canonical projection} of $W_0-((\partial {\bar{N}}_{0}) \times \{0\})$ onto ${\rm Int}\ {\bar{N}}_{0}$". 

We can define another set $$M_{0,t^{\prime}}:=\{(x,y) \in {\mathbb{R}}^n \times {\mathbb{R}}^{m-n+1}={\mathbb{R}}^{m+1} \mid F_{S^{m-n}}(y,t^{\prime})=0\}$$
 for $0<t^{\prime} \leq \frac{{\prod}_{j=1}^l (f_{{\partial}_j}(x))}{T}$. We present several properties based on our presented arguments.
 We can also know $M_{0,t^{\prime}}:=\{(x,y) \in \bar{N_0} \times {\mathbb{R}}^{m-n+1}={\mathbb{R}}^{m+1} \mid F_{S^{m-n}}(y,t^{\prime})=0\}$.
 For distinct real numbers $t^{\prime}=t_1, t_2$ here, the sets $M_{0,t_1}-((\partial \bar{N_0}) \times \{0\})$ and $M_{0,t_2}-((\partial \bar{N_0}) \times \{0\})$ are disjoint.
 These sets $M_{0,t^{\prime}}-((\partial \bar{N_0}) \times \{0\})$ are $m$-dimensional smooth manifolds with no boundaries. They are seen as the graphs of smooth functions as $M_0-((\partial \bar{N_0}) \times \{0\})$ and diffeomorphic to $M_0-((\partial \bar{N_0}) \times \{0\})$. Here some variable $y_j$ is regarded as the variable of the target of the function. We can also know that the manifold $M_{0,t^{\prime}}$ is an $m$-dimensional smooth and closed manifold and a non-singular real algebraic manifold. This manifold is also homeomorphic to the manifold $M_0:=M_{0, \frac{{\prod}_{j=1}^l (f_{{\partial}_j}(x))}{T}}$. We can show these facts by arguments similar to the corresponding arguments in the case of $M_{0}$ and $M_{0}-((\partial \bar{N_0}) \times \{0\})$. 
 
From our definitions, we also have the relations $W_0=({\bar{N}}_{0} \times \{0\}) \bigcup ({\bigcup}_{t^{\prime}  \in (0, \frac{{\prod}_{j=1}^l (f_{{\partial}_j}(x))}{T}]} M_{0,t^{\prime}})$ and $W_0=({\bar{N}}_{0} \times \{0\}) \sqcup ({\sqcup}_{t^{\prime}  \in (0, \frac{{\prod}_{j=1}^l (f_{{\partial}_j}(x))}{T}]} (M_{0,t^{\prime}}-((\partial {\bar{N}}_{0}) \times \{0\})))$.

We go back to the trivial bundle $W_0-((\partial {\bar{N}}_{0}) \times \{0\})$ over ${\rm Int}\ {\bar{N}}_{0}$ which is not seen as a smooth bundle. This is also regarded as a smooth bundle by so-called Ehresmann's fibration theorem (or a relative version of this). For this, remember that around a point of the form $(x,y) \in (({\rm Int}\ {\bar{N}}_{0}) \times {\mathbb{R}}^{m-n+1}) \bigcap M_0$, it is represented as the graph of a canonically obtained smooth function by regarding some variable $y_j$ as the variable of the target of the function. Remember also that the manifold $M_{0,t^{\prime}}$ enjoys a similar property, for example.

In the case $F_{S^{m-n}}(y,t)=t-||y||^2$, the previous smooth bundle $W_0-S(f_0)$ which may not be trivial in general is also a trivial smooth bundle. Our arguments here generalize ones from "STEP 3 in the proof of Main Theorems \ref{mthm:1} and \ref{mthm:2}".

These arguments also imply the relation $f_0 \circ{\Phi}^{\prime}=f^{\prime}$.
We can also discuss for the existence of the homeomorphisms of (\ref{mthm:3.6}) similarly. More precisely, to obtain a product-organized special generic map $f:M \rightarrow {\mathbb{R}}^n$ satisfying the relation $\bar{f}={\bar{f}}_{{\mathbb{R}}^n}$ and a homeomorphism $\phi:{\mathbb{R}}^n \rightarrow {\mathbb{R}}^n$ satisfying the relation $f^{\prime} =\phi \circ f$ and $\phi(f(S(f)))=f_0(S(f_0))$, we apply Proposition \ref{prop:3} with Remark \ref{rem:1} (\ref{rem:1.2}): remember also that the embeddings ${\bar{f}}_{{\mathbb{R}}^n}$ and ${\bar{f}}_{{\mathbb{R}}^n,0}$ are mutually isotopic smooth embeddings. 
In addition, the homeomorohism ${\Phi}^{\prime}$ is our desired homeomorphism $\Phi$.
This completes the proof of the statement similar to that of Main Theorem \ref{mthm:1} with the properties (\ref{mthm:3.5}, \ref{mthm:3.6}).

We can also check the property (\ref{mthm:3.7}) easily by our construction.
We can check that around the set $S(f_0)$, $f_0$ is represented as a natural product map as the case of "$f_0$ in Main Theorems \ref{mthm:1} and \ref{mthm:2}".

We can also check Main Theorem \ref{mthm:3} (\ref{mthm:3.8}) easily by our construction.
 A result close to this is also presented in the previous version (2nd version) of \cite{kitazawa5} (\cite[Theorem 1 (3)]{kitazawa5}). There "$\overline{D^{\prime}}$" is seen as the ($m+1$)-dimensional bounded and compact manifold $W_0$ bounded by $M_0$. The topology of the ($m+1$)-dimensional manifold is not so precisely discussed there. Only exposition on "homotopy" is presented there. Of course we do not assume arguments and result of \cite{kitazawa5}. 

This completes the proof.
\end{proof}

\begin{Ex}
	\label{ex:1}
The case $F_{S^{m-n}}(y,t)=t-||y||^2$ ($y \in {\mathbb{R}}^{m-n+1}$) gives one of simplest cases. For this we can choose $a$ as an arbitrary positive real number. 
The bundle $M_{{\rm B},0}$ in (\ref{mthm:3.1}) is regarded as a trivial linear bundle with ${\Delta}^k=D^k$ in this case. 
Furthermore, in this case, we can consider a suitable diffeomorphism ${\Phi}_{F_{{S}^{m-n}},a}$ and regard that (\ref{mthm:3.2}) can be discussed in the smooth category, for example. In short, we can argue all in the smooth situation in Main Theorem \ref{mthm:3} except arguments related to the deformation of the embedding ${\bar{f}}_{{\mathbb{R}}^n}$ of STEP 1. This is also for Main Theorems \ref{mthm:1} and \ref{mthm:2}. Main Theorems \ref{mthm:1} and \ref{mthm:2} give one of most important specific and simplest cases.

This is generalized for more general cases with $m-n \neq 4$ for example. This is thanks to fundamental arguments on singularity theory and differential topology. This is also thanks to well-known classical and sophisticated theory of diffeomorphisms on spheres. This is also related to studies on diffeomorphism groups of the unit spheres, studied in \cite{hatcher}, and so-called "pseudo-isotopy", pioneered in \cite{cerf}, for example. 

We also apply related arguments implicitly in some of our previous arguments in the present paper. See also Remark \ref{rem:5}, in the next section, for example. 

\end{Ex}

We give additional explicit cases as Examples. They are also a kind of our new remark.
\begin{Ex}
	\label{ex:2}
Consider the case of the form $F_{S^{m-n}}(y,t)=t-F_{S^{m-n},{\mathbb{R}}^{m-n+1}}(y)$ with a function $F_{S^{m-n},{\mathbb{R}}^{m-n+1}}$ satisfying the following. 
\begin{itemize}
	\item The restriction of $F_{S^{m-n},{\mathbb{R}}^{m-n+1}}$ to its zero set is a Morse function on an ($m-n+1$)-dimensional smooth manifold with no boundary. 

	\item The minimum of $F_{S^{m-n},{\mathbb{R}}^{m-n+1}}$ is $0$ and the condition ${F_{S^{m-n},{\mathbb{R}}^{m-n+1}}}^{-1}(0)=\{0\} \subset {\mathbb{R}}^{m-n+1}$ holds.
	\item The preimage $F_{S^{m-n},{\mathbb{R}}^{m-n+1}}(t)$ is a non-singular real algebraic manifold of dimension $m-n$ for $t>0$ and diffeomorphic to $S^{m-n}$.
\end{itemize}

Example \ref{ex:1} is generalized here.

We have a special generic map $f_0$ as in the case of Main Theorems \ref{mthm:1} and \ref{mthm:2}.
\end{Ex}
\begin{Ex}
	\label{ex:3}
We also note that Main Theorem \ref{mthm:3} respects the case $F_{S^{m-n},{\mathbb{R}}^{m-n+1}}(y):={\Sigma}_{j=1}^{m-n+1} k_{j,1} {y_j}^{2 k_{j,2}}$ with positive numbers $k_{j,1}$ and positive integers $k_{j,2}$, for example. See also \cite[Remark 3]{kitazawa2}.
Some examples of Example \ref{ex:2} are regarded as specific cases of the present cases.
\end{Ex}
\section{Additional comments on our new result.}

Remark \ref{rem:2} compares our main result to previous results by the author.
\begin{Rem}
\label{rem:2}
In \cite{kitazawa2, kitazawa3, kitazawa5}, we consider cases where
in Main Theorems, we do not need to change the embedding ${\bar{f}}_{{\mathbb{R}}^n}$ and can regard ${\bar{f}}_{{\mathbb{R}}^n}={\bar{f}}_{{\mathbb{R}}^n,0}$ with the polynomials $f_{{\partial}_j}$ being given beforehand. 
In short, only "STEP 2 in our proof of Main Theorems \ref{mthm:1} and \ref{mthm:2}" is considered and related short exposition on "STEP 3 in our proof of Main Theorems \ref{mthm:1} and \ref{mthm:2}" is presented there. The latter presentation is a kind of short remark. This does not affect main ingredients of \cite{kitazawa2, kitazawa3, kitazawa5}. 
"Our STEP 3" explicitly shows that $f_0$ is "topologically special generic" first.

These studies \cite{kitazawa2, kitazawa3, kitazawa5} also consider very explicit cases. For example, the polynomials are of degree $2$ and for spheres. These studies are also regarded as studies finding nice classes of real algebraic maps including the canonical projections of the unit spheres as simplest cases. 
\end{Rem}


The following is on another problem.

\begin{Rem}
\label{rem:3}
As related studies, we also present a comment on another problem, studying global structures of functions obtained by composing the canonical projections to the real algebraic maps in Main Theorems.
Some are closely related to Remark \ref{rem:2}, just before.

More precisely, such studies are on {\it Reeb graphs} of these functions. The {\it Reeb graph} of a smooth function is the graph obtained as the quotient space of the manifold of the domain and consists of all connected components of all preimages of single points. 
The vertex set of the Reeb graph consists of all elements representing connected components containing some singular points of the function. 

\cite{reeb} is a pioneering study considering Reeb graphs. Such tools are fundamental and strong tools in studying manifolds by using graphs. In 2006, \cite{sharko} has asked whether we can reconstruct a (nice) smooth function on a suitable manifold whose Reeb graph is isomorphic to a prescribed graph. 
\cite{masumotosaeki, michalak} also show explicit answers. \cite{kitazawa1} is motivated by such studies. There the author has succeeded in construction of nice smooth functions with prescribed preimages, first.
 
\cite{kitazawa2} considers this problem in the real algebraic category and constructs real algebraic functions with prescribed Reeb graphs in explicitly and nicely considered situations, first.
This is also regarded as a study finding nice classes of real algebraic functions including the canonical projection of the unit sphere into $\mathbb{R}$. For the Reeb graph of the canonical projection of the unit sphere of dimension at least $2$ into $\mathbb{R}$, its Reeb graph consists of exactly one edge and exactly two vertices. 
For this, \cite{bodinpopescupampusorea, sorea1, sorea2} are also closely related, help the author to do such studies, and motivate us to study further. They also present some important topics on singularity theory, real algebraic geometry and topology of non-singular real algebraic curves in ${\mathbb{R}}^2$ with regions surrounded by these curves in ${\mathbb{R}}^2$.
\end{Rem}

Remark \ref{rem:4} is on "STEP 1 in our proof of Main Theorems \ref{mthm:1} and \ref{mthm:2}". This is also independent of our present study.
Here we also give comments on singularity theory of smooth maps, Mather has pioneered and studied mainly in \cite{mather1, mather2, mather3, mather4, mather5, mather6, mather7, mather8}. See also \cite{ruas}. Remark \ref{rem:4} is also closely related to \cite{bodinpopescupampusorea}. This is on graphs in ${\mathbb{R}}^2$ which are generic from the viewpoint of singularity theory. This mainly studies regions surrounded by non-singular real algebraic curves and collapsing naturally to the graphs. 
We do not present elementary notions on graphs precisely. We except related elementary knowledge. We also note that we omit several rigorous discussions in Remark \ref{rem:4} and that for our main ingredients of the present paper, we do not need to understand this rigorously. We also give comments that this is also related to Remarks \ref{rem:2} and \ref{rem:3}.  
\begin{Rem}
\label{rem:4}
\cite{bodinpopescupampusorea} considers a finite graph $K$ enjoying the following properties.
\begin{itemize}
\item We can embed $K$ by a piecewise smooth embedding $e:K \rightarrow {\mathbb{R}}^2$.
\item The restriction of the function ${\pi}_{2,1} \circ e:K \rightarrow \mathbb{R}$ to each edge is a smooth embedding.
\item At distinct vertices, the values of ${\pi}_{2,1} \circ e$ are always distinct.
\item The degree of each vertex is $1$ or $3$.
\item At the vertex where the function ${\pi}_{2,1} \circ e:K \rightarrow \mathbb{R}$ has a local extremum, the degree of the vertex is $1$.
\end{itemize}
An {\it algebraic domain} is an open set in ${\mathbb{R}}^2$ surrounded by mutually disjoint non-singular connected real algebraic curves.
The {\it Poincar\'e-Reeb graph} of an algebraic domain $R$ is the graph defined as the quotient space of the closure $\overline{R}$ and the space consisting of all connected components of all preimages of all single points for the function ${\pi}_{2,1} {\mid}_{\overline{R}}$: the natural equivalence relation on the closure $\overline{R}$ is defined here. The vertex set of the graph can be and is defined as the set of all points for connected components containing some singular points of the function ${\pi}_{2,1} {\mid}_{\partial \overline{R}}$.
The paper says that for the embedded graph $e(K)$, we can have an algebraic domain $R_K$ enjoying the following.
\begin{enumerate}
\item $e(K) \subset R_K$.
\item The Poincar\'e-Reeb graph of $R_K$ is isomorphic to and naturally regarded as $e(K)$. 
\item We also have the situation in the following. Of course $\overline{R_K}$ is the closure of $R_K$ in ${\mathbb{R}}^2$.
\begin{enumerate}
\item The restriction of ${\pi}_{2,1}$ or ${\pi}_{2,1} {\mid}_{\overline{R_K}}:\overline{R_K} \rightarrow \mathbb{R}$ to the boundary $\partial \overline{R_K}$ is a Morse function such that at distinct singular points the values are always distinct.
\item $\partial \overline{R_K}$ is regarded as the intersection of some bounded open set and the zero set of some real polynomial function. 
\end{enumerate}
\end{enumerate}
Our desired region $R_K$ can be also obtained by a kind of approximation of the Weierstrass type in the class $C^2$.
According to arguments similar to ones in \cite[Special Case 5]{kollar}, we can consider the $C^2$ approximation here in such a way that our $\partial \overline{R_K}$ here is regarded as the zero set of another real polynomial function.
 
\end{Rem}
The following is on Main Theorem \ref{mthm:3} and related to some theory of diffeomorphism groups and differential topolgy of manifolds.
\begin{Rem}
\label{rem:5}
	In Main Theorem \ref{mthm:3}, let $m=n$. By respecting the well-known facts that the diffeomorphism group of $D^1$ has the homotopy type of the two-point discrete set and that each point of the two-point set corresponds to orientation preserving diffeomorphisms and orientation reversing ones, respectively, $M$ and $M_0$ are shown to be diffeomorphic to the so-called double of the manifold $\bar{N_0}$.
	\end{Rem}

The following is closely related to one of central problems in real algebraic geometry. This may be closely related to our future work.
\begin{Rem}
\label{rem:6}
Construction of real algebraic curves in the real affine space ${\mathbb{R}}^2$ and real algebraic hypersurfaces of a general dimension has been difficult and important in real algebraic geometry. 
This comes from Hilbert's 16th problem and for example, Viro's sophisticated theory such as \cite{gabard, viro} presents strong tools.
At present, even explicit construction and classifications for cases of real algebraic curves are difficult. 
Our pioneering studies of the construction respect the unit spheres and the canonical projections of them and we do not restrict attention to the $1$-dimensional case. Our pioneering studies are also new in concentrating not only on the hypersurfaces themselves but also on their canonical projections.
\end{Rem}
The author recognizes that our construction shows explicit cases of Viro's methods. In short, Viro's methods glue local hypersurfaces together. We can have non-singular real algebraic hypersurfaces by deforming and smoothing (Viro's patchworking). In our paper we do not need to understand this theory well. However, to study more, we need to understand the theory more.

\section{Acknowledgement.}
This new study is based on studies on special generic maps and explicit construction of real algebraic functions and maps of the author.
One of related events motivating the author to present this study is the conference "Extension of the singularity theory" (https://www2.akita-nct.ac.jp/kasedou/workshop/rims2023/index.html)" and related discussions on singularity theory, differential geometry, and real algebraic geometry.
See also the previous versions \cite{kitazawa6}.

\end{document}